%% file: main.tex
\documentclass{IEEEtran4PSCC}
\ifCLASSINFOpdf
   \usepackage[pdftex]{graphicx}
\else
   \usepackage[dvips]{graphicx}
\fi

\usepackage[cmex10]{amsmath}
\usepackage{float}
\floatstyle{ruled}
\newfloat{model}{thp}{lop}
\floatname{model}{Model}

\usepackage{url}
\usepackage{booktabs}
\usepackage{amsmath}
\usepackage{amssymb}
\usepackage{bm}
\usepackage{xcolor}
\usepackage{tablefootnote}
\usepackage[colorlinks=false, linkbordercolor=red, citebordercolor=green, urlbordercolor=white]{hyperref}
\usepackage{multirow}
\usepackage{algorithm}
\usepackage[noend]{algpseudocode}

\input{tex/definitions}

\usepackage{todonotes}

\makeatletter
\let\old@ps@headings\ps@headings
\let\old@ps@IEEEtitlepagestyle\ps@IEEEtitlepagestyle
\def\psccfooter#1{%
    \def\ps@headings{%
        \old@ps@headings%
        \def\@oddfoot{\strut\hfill#1\hfill\strut}%
        \def\@evenfoot{\strut\hfill#1\hfill\strut}%
    }%
    \def\ps@IEEEtitlepagestyle{%
        \old@ps@IEEEtitlepagestyle%
        \def\@oddfoot{\strut\hfill#1\hfill\strut}%
        \def\@evenfoot{\strut\hfill#1\hfill\strut}%
    }%
    \ps@headings%
}
\makeatother

\usepackage{hyperref}

\begin{document}
\title{Volt/VAR Optimization in Transmission Networks with Discrete-Control Devices}

\author{\IEEEauthorblockN{Shuaicheng Tong\IEEEauthorrefmark{1},
Michael A. Boateng\IEEEauthorrefmark{2},
Mathieu Tanneau\IEEEauthorrefmark{1},
Pascal Van Hentenryck\IEEEauthorrefmark{1}}
\IEEEauthorblockA{\IEEEauthorrefmark{1}School of Industrial and Systems Engineering, Georgia Institute of Technology}
\IEEEauthorblockA{\IEEEauthorrefmark{2}School of Electrical and Computer Engineering, Georgia Institute of Technology}

\{stong38, mboateng6\}@gatech.edu, \{mathieu.tanneau,pascal.vanhentenryck\}@isye.gatech.edu}

\maketitle


\begin{abstract}
Voltage (Volt) and reactive-power (VAR) control in transmission networks is critical for reliability and increasingly needs fast, implementable decisions. This paper presents a transmission Volt/VAR Optimization (VVO) framework that co-optimizes discrete control of on-load tap-changing transformers (OLTCs) and capacitor banks (CBs) with AC power flow (ACPF) physics to improve voltage stability and minimize VAR generation.
The framework follows a relax–round–resolve pipeline: a continuous relaxation proposes targets, a rounding step selects feasible discrete settings, and a final solve enforces AC power flow physics. Extensive experiments on IEEE, PEGASE, and RTE systems show consistent improvements in voltage and VAR quality metrics with modest generator redispatch while preserving economic operation and achieving compatible runtimes with real-time transmission operations.
 \end{abstract}


% Use this to place sponsorships
% \thanksto{\noindent Submitted to the 24th Power Systems Computation Conference (PSCC 2026).}

% \section*{Nomenclature}

% \subsection{Abbreviations}
% \begin{tabbing}
%     OPF \hspace{1.5cm} \= Optimal Power Flow \\
%     VVO \> Volt/VAR Optimization \\
%     VAR \> Volt-Ampere Reactive \\
%     DER \> Distributed Energy Resource \\
%     MINLP \> Mixed Integer Nonlinear Programming \\
%     CB \> Capacitor Bank \\
%     OLTC \> On-Load Tap Changer Transformers \\
% \end{tabbing}

\input{tex/introduction}

\input{tex/formulation}

\input{tex/methodology}

\input{tex/experiments}

\input{tex/conclusion}

% The 'triggered' command can be changed if desired:

\bibliographystyle{IEEEtran}
\bibliography{refs}

\clearpage

\end{document}

%% file: tex/introduction.tex
\section{Introduction}

\subsection{Motivation}
Power system operators must keep voltages within tight bounds while transferring large amounts of real power across transmission networks. \cite{Stott2012OptimalPF}
Voltage control is essential for equipment protection, power quality, and system stability as deviations outside acceptable ranges can lead to cascading voltage collapse \cite{BickfordVARControl, Kundur2004_StabilityDefs}.
Achieving this in real time requires effective Volt/VAR control: sufficient local VAR support stabilizes voltage profiles, limits reactive power circulation and losses, as well as preserves thermal headroom for active-power transfers \cite{JAY2021110967}.

In modern grids, Transmission System Operators (TSOs) coordinate voltage control with power dispatch to maintain secure and economic operations.
This coordination relies on discrete devices such as on-load tap-changing transformers (OLTCs), capacitor banks (CBs), and, where available, flexible AC Transmission System (FACTS) devices \cite{CAISO2019,IEEE1547_2018, Kundur1994PowerSS}.
Because OLTCs and CBs act in discrete increments within finite limits, coordinating them under AC power-flow physics is combinatorially challenging in computation \cite{BienstockVerma2015, CAPITANESCU201657}.
As a result, operators rely on rule-based adjustments, contingency studies, and sensitivity analyses \cite{ElizondoSurvey2017_TransVoltageControl} that can under-utilize device flexibility or over-correct in stressed conditions. It is hence essential to strengthen the coordination between power dispatch and voltage/reactive-power control so that operational decisions remain consistent with network physics.

\subsection{Literature Review}
Voltage and reactive-power management are intertwined with network physics. It shapes transfer capability (maximum AC-feasible real-power transfer subject to thermal, voltage, and stability limits), reliability, and operating economy \cite{Kundur1994PowerSS}. When generators or devices reach VAR limits, stability margins shrink and feasible transfers contract \cite{Molzahn2013_PowerFlowInsolvability,Zhao2008_BusTypeSwitching}.
Computationally, the VVO problem is typically posed as a mixed-integer nonlinear program (MINLP) \cite{Momoh1997_ChallengesOPF,PapalexopoulosLargeOPF}. Prior approaches of solving MINLPs in power system include rounding schemes \cite{Macfie2010ShuntRounding,Karoui2008_EEM_IPOPF}, relaxations and iterative cutting plane methods \cite{Molzahn2019Relaxation,Aigner2023ACOPFGlobalOptimality,PlatbroodGenericApproaches, Ding2004_IMA_IPCuttingPlane_OPF}, homotopy/continuation methods \cite{agarwal2018continuous,McNamara2022TwoStageHomotopy}, genetic algorithm/population-based heuristics \cite{mataifa_vvo_survey, Chamba2023_Energies_ORPD_MAS} and primal-dual interior-point methods \cite{Liu2002ExtendedIP}.
However, a gap remains: an \emph{AC-feasible, transmission-scale} VVO formulation that integrates operation and realistic device setpoints. Much prior work targets distribution systems (often without fixed operating setpoints and at smaller scales) \cite{Jha2020VoltControl,ortmann2020experimental}. Linearizations of MINLP into mixed-integer linear program (MILP) struggle when AC nonlinearities interact with discrete OLTC/CB controls, yielding infeasibilities or poor performance at scale \cite{Aigner2023ACOPFGlobalOptimality, McNamara2022TwoStageHomotopy}. Therefore, there is a need for methods that preserve AC physics, respect discrete device settings, and scale within time budgets.

\subsection{Contribution}

To address these challenges, the paper formulates the VVO problem as a MINLP and solves it via a “relax–round–resolve” framework, which is well-established \cite{Arora1994ReviewDiscrete, TinneyDeficienciesInOPF}. Simultaneously, the pitfalls of a naïve relax-and-round are well known: \emph{(i)} rounding can yield sub-optimal device schedules as the relaxed stage objective is unaware of the eventual discretization \emph{(ii)} many studies assume unrealistic VAR-control step sizes or initial values, which leads to settings that are hard to implement in practice \cite{McNamara2022TwoStageHomotopy,Capitanescu2010SensitivityOPF}. The paper designs around these issues by modeling constraints in common manufacturer-defined operating ranges, rather than idealized controls.

The paper's contributions are summarized as follows:
\begin{itemize}
  \item \textit{An operations-aligned objective} that maintains voltage and reactive-power setpoints while penalizing active-power redispatch, yielding minimal departures from reference setpoints.
  \item \textit{A realistic device model} with discrete OLTC tap ratios and shunt CB modules reflecting manufacturer step sizes and operating limits in transmission systems, including standardized initial taps to ensure implementable starting points.
  \item \textit{A scalable relax-round-resolve workflow} that returns AC-feasible solutions and operator-ready discrete settings across diverse scales of transmission systems.
\end{itemize}

%% file: tex/formulation.tex
\section{Problem Formulation}
\label{sec:prob_form}

In all that follows, consider a power grid represented as a simple directed graph $G \, {=} \, (\NODES, \EDGES)$, where $\NODES \, {=} \, \{1, ..., N\}$ is the set of buses, and $\EDGES \, {\subseteq} \, \NODES {\times} \NODES$ is the set of branches, i.e., lines and transformers.
The set of branches leaving (resp. entering) bus $i$ is denoted by $\EDGES^{+}_{i}$ (resp. $\EDGES^{-}_{i}$).
The set of generators is denoted $\GENERATORS$.
For ease of reading, assume that each bus has exactly one generator and one load, and that there are no parallel branches.
This is purely a notational assumption: the proposed methodology and the paper's numerical experiments support parallel branches and an arbitrary number of devices at each bus.

The imaginary unit is denoted by $\im \, {\in} \, \mathbb{C}$, i.e., $\im^{2} \, {=} \, {-}1$.
The complex conjugate of $z \, {=} \, x \, {+} \, \im y , {\in} \, \mathbb{C}$ is $z^{\star} \, {=} \, x \, {-} \, \im y$.
Unless specified otherwise, and with the exception of $\im$, bold (resp. non-bold) symbols denote decision variables (resp. input parameters).
The complex voltage, generation and demand at bus $i \, {\in} \, \NODES$ are denoted by $\V_{i} \, {=} \, \VM_{i} \angle \VA_{i}$, $\SG_{i} \, {=} \, \PG_{i} \, {+} \, \im \QG_{i}$ and $\Sd_{i} \, {=} \, \pd_{i} \, {+} \, \qd_{i}$, respectively.
The forward and reverse complex power flows on branch $ij \, {\in} \, \EDGES$ are denoted by $\SF_{ij} \, {=} \, \PF_{ij} \, {+} \, \im \QF_{ij}$ and $\ST_{ij} \, {=} \, \PT_{ij} \, {+} \, \im \QT_{ij}$, respectively.

\subsection{Capacitor Bank Modeling}
\label{subsec:CB_modeling}
In transmission networks, shunt elements represent current injection to or absorption from the ground at a bus. 
Capacitor banks (CBs) allow operators to control nodal VAR injections by adjusting the shunt admittance.
Each CB comprises several modules that can be switched on or off to adjust the device susceptance in discrete increments, thereby increasing reactive currents, altering reactive power flows, and regulating bus voltages.

Whereas shunt elements are modeled as fixed admittances in canonical ACOPF formulations, this paper explicitly models CB operations.
Namely, the total shunt admittance at bus $i$ is
\begin{align}
\label{eq:ys_nodal_admittance}
\ys_i(\CB_i) = \gs_i + \im(\bs_i + \CB_i)
\end{align}
where $\gs_{i} + \im \bs_{i}$ denotes the shunt element's reference admittance, and $\CB_{i}$ denotes the additional susceptance contributed by the active CB modules. Shunt conductance, $\gs_i$, typically represents fixed losses to the ground, which is not controllable in practice and hence not considered in this study.
$\CB_i$ takes values in the discrete set $\CBNODES_{i}$, corresponding to the device's available setpoints.
If bus $i$ is not equipped with CBs, then $\CBNODES_i \, {=} \, \{0\}$.
Accordingly, Kirchhoff’s current law at bus $i$ is formulated, in complex form, as
\begin{align}
\label{eq:kcl}
    \SG_i
    - \Sd_i
    - (\ys_i(\CB_i))^{\star} |\V_i|^{2}
    = 
    \sum_{e\in\Eiplus}\SF_{e}
    + \sum_{e\in\Eiminus}\ST_{e}.
\end{align}

\subsection{Transformer Modeling}
\label{subsec:trafo_modeling}
On-load tap-changing transformers allow system operators to adjust network voltages as loading conditions evolve. By adjusting the tap ratio, operators tune the electrical coupling between connected buses to meet operational objectives.

Transformers provide voltage control between buses by scaling voltage magnitudes and shifting voltage phases. Operators change transformer tap ratios to keep nodal voltages within limits and reduce losses. This also modifies reactive and active power flows without large-scale generator redispatch. In steady-state power-flow analysis, the transformer’s effect is represented by the $2\times2$ branch admittance submatrix between the from-bus $i$ and to-bus $j$.

\begin{align}
\label{eq:Ybr_matrix}
\Ybr_{ij}(\TAP_{ij}) =
\begin{bmatrix}
    {\yff_{ij}} \times {\TAP_{ij}^{-2}}
        & {\yft_{ij}} \times {\TAP_{ij}^{-1}} \\
    {\ytf_{ij}} \times {\TAP_{ij}^{-1}}
        & \ytt_{ij}
\end{bmatrix}
\in \mathbb{C}^{2 \times 2},
\end{align}
where $\TAP_{ij} \, {>} \, 0$ denotes the branch's tap ratio, which takes values in $\mathcal{T}_{ij} \, {\subseteq} \, \mathbb{R}$.
If branch $ij$ is a line, then $\mathcal{T}_{ij} \, {=} \, \{1\}$.
For a transformer, $\mathcal{T}_{ij}$ is typically a discrete set, corresponding to discrete setpoints of the underlying device.
Branch phase-shift angles are assumed to be fixed, and their contribution to the branch admittance is therefore captured by $\yft_{ij}, \ytf_{ij}$.

Ohm's law for branch $ij \, {\in} \, \EDGES$ is then stated as
\begin{align}
\label{eq:ohm:complex}
    \begin{bmatrix}
        \SF_{ij}\\
        \ST_{ij}
    \end{bmatrix}
    &=
    \begin{bmatrix}
        \V_{i} & \\
        & \V_{j}
    \end{bmatrix}
    \times
    \left(\Ybr_{ij}(\TAP_{ij})\right)^{\star}
    \times
    \begin{bmatrix}
        \V_{i}^{\star}\\
        \V_{j}^{\star}
    \end{bmatrix}
    ,
\end{align}
and thermal constraints read
\begin{align}
    \label{eq:thermal}
    |\SF_{ij}|, |\ST_{ij}| &\leq \smax_{ij},
\end{align}
where $\smax_{ij}$ is the branch's maximum apparent power flow.

\iffalse{
\begin{subequations}
\label{eq:ohm}
\begin{align}
    \label{eq:ohm:pf}
    \PF_{ij} &=
        \gij \frac{\VM_i^{2}}{\TAP_{ij}^{2}}
        - \gij\frac{\VM_i \VM_j}{\TAP_{ij}}\cos\VA_{ij}
        - \bij\frac{\VM_i \VM_j}{\TAP_{ij}}\sin\VA_{ij}
        \\
    \label{eq:ohm:qf}
    \QF_{ij} & =
        - (\bij+\frac{\bcij}{2})\frac{\VM_i^{2}}{\TAP_{ij}^{2}}
        + \bij\frac{\VM_i \VM_j}{\TAP_{ij}}\cos\VA_{ij}
        - \gij\frac{\VM_i \VM_j}{\TAP_{ij}}\sin\VA_{ij}
        \\
    \label{eq:ohm:pt}
    \PT_{ij}  & =
        \gij\VM_j^{2}
        - \gij\frac{\VM_i \VM_j}{\TAP_{ij}}\cos\VA_{ij}
        + \bij\frac{\VM_i \VM_j}{\TAP_{ij}}\sin\VA_{ij}
    \\
    \label{eq:ohm:qt}
    \QT_{ij} & =
        - (\bij+\frac{\bcij}{2})\VM_j^{2}
        + \bij\frac{\VM_i \VM_j}{\TAP_{ij}}\cos\VA_{ij}
        + \gij\frac{\VM_i \VM_j}{\TAP_{ij}}\sin\VA_{ij}
\end{align}
\end{subequations}
where the branch angle difference $\VA_{ij} \, {=} \, \VA_{i} \, {-} \, \VA_{j}$.\\
}\fi

\subsection{Volt-Var Optimization Model}
\label{subsec:vvo_model}

\begin{model}[!t]
    \caption{Volt-Var Optimization Model}
    \label{model:VVO}
    \begin{subequations}
    \label{eq:vvo}
    \begin{align}
        \min_{\V, \SG, \SF, \ST, \CB, \TAP} \quad
        & \label{eq:vvo_obj}
            % \lamv \psi_{v}(\VM) + \lamq  \psi_{q}(\QG) + \psi_{p}(\PG)\\
            \psi(\VM, \QG, \PG)\\
        \text{s.t.} \quad
        & \eqref{eq:kcl}, \eqref{eq:ohm:complex}, \eqref{eq:thermal} && \nonumber \\
        &
          \thetaref = 0
        &&
        \label{constraint:slack_bus_angle} \\
        &
          \pgmin \le \PG_{i} \le \pgmax
        && \forall i\in\GENERATORS
        \label{constraint:nodal_active_gen_limit} \\
        &
          \qgmin \le \QG_{i} \le \qgmax
        && \forall i\in\GENERATORS
        \label{constraint:nodal_reactive_gen_limit} \\
        &
          \dvamin_{ij} \le \VA_{i} - \VA_{j} \le \dvamax_{ij}
        && \forall ij \in \EDGES
        \label{constraint:br_angle_deviation} \\
        &
          \vmmin \le \VM_{i} \le \vmmax
        && \forall i\in\NODES
        \label{constraint:nodal_vm_limit} \\
        &
          \CB_{i} \in \CBNODES_i
        && \forall i\in\NODES
        \label{constraint:CB_modules} \\
        &
          \TAP_{ij} \in \TAPEDGES_{ij}
        && \forall ij\in\EDGES
    \label{constraint:tap_position}
    \end{align}
    \end{subequations}
\end{model}

The proposed Volt-Var-Optimization (VVO) is formulated as a Mixed-Integer Nonlinear Programming (MINLP) problem in Model \ref{model:VVO}.
Constraint \eqref{constraint:slack_bus_angle} sets the voltage angle of the slack bus to zero.
Constraints \eqref{constraint:nodal_active_gen_limit} and \eqref{constraint:nodal_reactive_gen_limit} impose minimum and maximum limit on active and reactive power generation, respectively.
Constraint \eqref{constraint:br_angle_deviation} limits the phase angle difference between a branch's two endpoints.
Constraint \eqref{constraint:nodal_vm_limit} ensures that voltage magnitude remain within an acceptable range.
Finally, constraints \eqref{constraint:CB_modules} and \eqref{constraint:tap_position} represent the discrete operating setpoints of capacitor banks and transformer taps, respectively.

The paper considers an objective function \eqref{eq:vvo_obj} of the form
\begin{align*}
    \psi(\VM, \QG, \PG)
    &=
    \lamv \psi_{v}(\VM)
    + \lamq \psi_{q}(\QG)
    + \lamp \psi_{p}(\PG)
    + \lamc \psi_{c}(\PG),
\end{align*}
where $\lamv, \lamq, \lamp, \lamc$ are non-negative weight parameters and functions $\psi_{v}, \psi_{q}, \psi_{p}, \psi_{c}$ are defined as
\begin{subequations}
\label{eq:vvo_obj_terms}
\begin{align}
  \label{eq:obj:v}
  \psi_{v}(\VM) 
    &= \sum_{i\in\NODES} (\VM_i - \vref_{i})^2, \\
  \label{eq:obj:q}
  \psi_{q}(\QG) 
    &= \sum_{i\in\NODES} (\QG_i - \qgref_{i})^2, \\
  \label{eq:obj:pref}
  \psi_{p}(\PG) 
    &= \sum_{i\in\NODES} (\PG_i - \pgref_{i})^2,
  \\ 
  \label{eq:obj:pcost}
  \psi_{c}(\PG)
    &= \sum_{i\in\NODES} \cost_{i}(\PG_{i}).
\end{align}
\end{subequations}
Thereby, $\psi_{v}, \psi_{q}$ and $\psi_{p}$ in \eqref{eq:obj:v}-\eqref{eq:obj:pref} penalize deviations from reference setpoints for voltage magnitude, reactive power dispatch and active power dispatch, respectively.
The term $\psi_{c}$ in \eqref{eq:obj:pcost} captures the economic cost of active power generation.
Unless specified otherwise, for every bus $i$, the paper considers $\vref_{i} = 1$, i.e., voltage magnitudes should be close to 1p.u., and $\qgref_{i} = 0$, i.e., reactive generation should be as small as possible.
The reference active power dispatch $\pgref$ is chosen as the active power dispatch set by the market.
Note that one can emulate a constraint of the form $\PG = \pgref$ by setting $\lamp$ to an infinite value.
Also note that, by setting $\lamv = \lamq =\lamp = 0$, and fixing variables $\CB$ and $\TAP$ to a reference value, Model \ref{model:VVO} reduces to the canonical ACOPF problem.

%% file: tex/methodology.tex
\section{Solution Methodology}
\label{sec:methodology}

Model \ref{model:VVO} can be solved by off-the-shelf MINLP solvers such as BARON or Gurobi \cite{gurobi}.
Nevertheless, during preliminary experiments on small and medium power grids, Gurobi failed to produce any feasible solution within an hour of computing time, thereby illustrating the numerical and combinatorial challenge of solving Model \ref{model:VVO} directly.
This is likely caused by numerical challenges, poor relaxations, and the problem's combinatorial structure brought by discrete variables $\CB$ and $\TAP$.
This behavior was observed even when restricting the range of discrete devices, thus suggesting that the non-convex nature of power flow constraints is particularly challenging.

\begin{algorithm}[!t]
\caption{VVO Solution Heuristic}
\label{alg:vvo_round_resolve}
\begin{algorithmic}[1]
\State \textbf{Input:} Network $(\NODES,\EDGES)$; ACOPF solution $\hat{\PFSOL}$; device settings $\CBNODES_i,\TAPEDGES_{ij}$
\Statex
\State \emph{Relaxed VVO}
\State Relax constraints \eqref{constraint:CB_modules} and \eqref{constraint:tap_position} to continuous intervals
\State Solve Model \ref{model:VVO} with relaxed constraints
\State Retrieve fractional solution $\TAPrvvo,\CBrvvo$
\Statex
\State \emph{Round device settings}
\For{each bus $i \in \NODES$}
  \State $\overline{\CB}_i \gets \mathrm{round}(\CBrvvo_i)$
\EndFor
\For{each branch $ij \in \EDGES$}
  \State $\overline{\TAP}_{ij} \gets \mathrm{round}(\TAPrvvo_{ij})$
\EndFor
\Statex
\State \emph{Fixed-device VVO}
\State Fix $\TAP_{ij} \equiv \overline{\TAP}_{ij}$ and $\CB_i \equiv \overline{\CB}_i$ in Model \ref{model:VVO}
\State Resolve Model \ref{model:VVO}
\If{feasible}
    \State \Return $(\overline{\TAP},\overline{\CB})$
\Else
    \State \Return no solution found
\EndIf
\end{algorithmic}
\end{algorithm}

To address this limitation the paper proposes a \emph{relax–round–resolve} ($R^{3}$) heuristic, presented in Algorithm \ref{alg:vvo_round_resolve}.
First (\emph{relax}), the continuous relaxation of \eqref{eq:vvo} is solved as a nonlinear non-convex continuous problem, yielding a fractional solution $\tilde{\CB}, \tilde{\TAP}$.
This step was found to be tractable using an off-the-shelf nonlinear optimization solver like Ipopt \cite{ipopt}.
Second (\emph{round}), the fractional solution is rounded to the nearest integer-feasible value, yielding a candidate solution $\bar{\CB}, \bar{\TAP}$ that satisfies constraints \eqref{constraint:CB_modules}-\eqref{constraint:tap_position}.
Third (\emph{resolve}), variables $\CB, \TAP$ in Model \ref{model:VVO} are fixed to $\bar{\CB}, \bar{\TAP}$, and the resulting problem is resolved.
Note that this last problem no longer contains any discrete variables, and is therefore amenable to a nonlinear optimization solver like Ipopt.
If this resolve fails to converge to a feasible solution, the heuristic terminates with a ``no solution found" code.

This heuristic framework can be augmented with a local search algorithm to further improve solution quality or remediate possible infeasibilities.
Nevertheless, as demonstrated in the numerical experiments of Section \ref{sec:num_exp}, the proposed heuristic finds feasible solutions most of the time, even on large-scale networks.

%% file: tex/experiments.tex
\section{Numerical Experiments}
\label{sec:num_exp}
\newcommand{\ieeeXXS}{\texttt{14\_ieee}}
\newcommand{\ieeeXS}{\texttt{30\_ieee}}
\newcommand{\ieeeS}{\texttt{57\_ieee}}
\newcommand{\ieeeM}{\texttt{118\_ieee}}
\newcommand{\ieeeL}{\texttt{300\_ieee}}

\newcommand{\pegaseXS}{\texttt{1354\_pegase}}   % 1.3k
\newcommand{\pegaseS}{\texttt{2869\_pegase}}    % 2.8k
\newcommand{\pegaseM}{\texttt{8387\_pegase}}    % 8.3k
\newcommand{\pegaseL}{\texttt{9241\_pegase}}    % 9.2k
\newcommand{\pegaseXL}{\texttt{13659\_pegase}}  % 13.7k

\newcommand{\rteXS}{\texttt{1888\_rte}}     % 1.9k
\newcommand{\rteS}{\texttt{2848\_rte}}     % 2.8k
\newcommand{\rteM}{\texttt{6468\_rte}}     % 6.5k variant A
\newcommand{\rteL}{\texttt{6470\_rte}}    % 6.5k variant B
\newcommand{\rteXL}{\texttt{6495\_rte}}   % 6.5k variant C
\newcommand{\rteXXL}{\texttt{6515\_rte}}  % 6.5k variant D

Unless specific otherwise, all optimization models are implemented in Julia 1.11.3 using JuMP \cite{JuMP} and solved with Ipopt \cite{ipopt}.
Experiments are conducted on dual Intel Xeon 6226@2.7GHz machine running Linux on the Phoenix cluster \cite{PACE}.
Each job is allocated 8 CPU cores, 32GB of RAM, and a 4-hour wallclock time limit.

A separate set of experiments were conducted where Model \ref{eq:vvo} was solved as a non-convex MINLP using Gurobi 12.0.3 \cite{gurobi}.
These experiments are omitted for sake of space, because Gurobi was unable to find a feasible solution in all tested cases.

\subsection{Experiment Setup}
\label{subsec:constraint_modeling}

The paper conducts experiments on small-to-large cases from the \texttt{PGLib} library \cite{pglib}, ranging from 118 to 13659 buses.
The \texttt{rte} and \texttt{pegase} cases were selected as they represent real-life transmission grids from the European system.
Table \ref{tab:devices} reports relevant case statistics, namely: the number of buses (\#bus), generators (\#gen) and capacitor banks (\#CB) in the grid, and the number of transmission lines (\#line) and transformers (\#tsfm).

\begin{table}[!t]
    \centering
    \caption{PGLib case statistics}
    \label{tab:devices}
    \input{tables/device_composition}
\end{table}

Capacitor bank and adjustable transformer taps are modeled as follows.
Buses where a shunt element is attached are assumed to be equipped with a capacitor bank.
The paper assumes that each capacitor bank consists of three modules, each with a rating of 10MVAr at nominal voltage.
Noting that all PGLib cases use a per-unit convention of $100$MVA, the paper thus uses $\CBNODES_{i} = \{0, 0.1, 0.2, 0.3\}$ for buses equipped with a CB, and $\CBNODES_{i} = \{0\}$ otherwise.
It is further assumed that each non-zero bus shunt admittance in the case data corresponds to one module of a capacitor bank being active.
This assumption is made to evaluate the performance of the proposed methodology on cases where VAR injection can be modified upwards and downwards.
Each transformer is assumed to have an adjustable tap ratio, with 33 possible tap positions $\{-16, ..., +16\}$, where the neutral position $0$ corresponds to a ratio $\tau = 1$ and each step corresponds to a $0.625\%$ increment.
This yields $\TAPEDGES_{e} = \{0.9, 0.90625, ..., 1.1\}$ for each transformer $e$.

Experiments are thus conducted as follows.
First, the PGLib casefile is loaded and pre-processed using the \texttt{make\_basic\_network} function in \texttt{PowerModels.jl} package.
Second, transformer tap ratio values are rounded to correspond to a valid tap setpoint, if necessary.
Third, an AC Optimal Power Flow is solved to obtain an AC-feasible reference setpoint.
Note that, in practice, such a setpoint may be the output of a market clearing step, a previous operating point, or based on a forecast of future operating conditions.
Fourth, the proposed heuristic is executed to find capacitor bank and transformer tap ratio setpoints that improve the system's volt-var profile.

All experiments use the objective function described in Section \ref{subsec:vvo_model}, with $\lambda_{v} \, {=} \, \lambda_{q} \, {=} \, 1$.
The paper also sets $\lambda_{c} \, {=} \, 1$, as this was found to improve convergence and prevent unacceptable increases in cost when active generation setpoints can be adjusted.
Noting that TSO often avoid active power re-dispatch for operational and economic reasons, the paper conducts a set of experiments where active power dispatches are fixed to their value in the reference setpoint, with the exception of the slack generator.
This is implemented by adding the constraint $\PG_{i} \, {=} \, \pgref_{i}$ to Model \ref{model:VVO} for all buses $i$ but the slack bus.
For ease of reading, the corresponding results are indicated by $\lambda_{p} = \infty$.
Finally, the paper conducts experiments with varying ranges of device setpoints, representing operational realities where large device movements can induce stability issues and accelerate equipment wear and tear.
Thereby, three ranges are considered: i) tap ratios can deviate by at most $\pm 3$ steps from their reference setpoint, and up to 2 CB modules may be activated, ii) tap ratio deviation of at most $\pm 3$ from the reference setpoint, and all (3) CB modules may be activated, and iii) full $\pm 16$ range of tap ratio and all CB modules may be activated.
This last setting is less realistic, and provides insight into the computational performance of the propose heuristic when dealing with a large combinatorial space.

\subsection{Performance Metrics}
\label{subsec:eval_metrics}

The paper uses the following metrics to evaluate the quality of solutions.
The quality of voltage profiles is measured via
\begin{align}
  \label{eq:vmae}
  \text{MAE}_v
  &= \frac{1}{|\NODES|} \sum_{i\in\NODES} \left| \VM_i - 1 \right|,
\end{align}
which reflect that voltage levels should be close to 1p.u.
The quality of reactive power dispatches is measured by
\begin{align}
    \label{eq:qmae}
    \text{MAE}_q
        &= \frac{1}{|\GENERATORS|}
            {\sum_{i\in\GENERATORS} \left| \QG_i \right|}
    ,
\end{align}
which capture the fact that reactive generation should be as close to zero as possible.
The paper also evaluates how solutions deviate from reference active generation setpoints, quantified by
\begin{align}
    \label{eq:pmae}
    \Delta\PG
    &=
    \frac{1}{|\GENERATORS|} {\sum_{i\in\GENERATORS} \left| \PG_i - \pgref_{i} \right|},
\end{align}
and the economic cost of adjusting active power dispatches
\begin{align}
    \%\Delta \cost
  &= 100 \cdot \frac{\psi_{c}(\PG)-\psi_{c}(p_{\mathrm{ref}}^{\mathrm{g}})}{\psi_{c}(p_{\mathrm{ref}}^{\mathrm{g}})}
\end{align}
where $\psi_{c}(\PG)$ denotes the cost of generation.
The heuristic's computational performance is measured by reporting the runtime of the relaxed and fixed steps of Algorithm \ref{alg:vvo_round_resolve}.
These computing times do not include the time for the initial ACOPF solve, as this is assumed to be a sank cost.

\begin{table}[!t]
    \centering
    \caption{Solution Quality of Pipelines of IEEE 118}
    \label{tab:118}
    \resizebox{\columnwidth}{!}{%
        \input{tables/results_118}
    }
\end{table}

\begin{table}[!t]
    \centering
    \caption{Solution Quality of Pipelines of PEGASE 1354}
    \label{tab:1354}
    \resizebox{\columnwidth}{!}{%
        \input{tables/results_1354}
    }
\end{table}

\begin{table}[!t]
    \centering
    \caption{Solution Quality of Pipelines of RTE 1888}
    \label{tab:1888}
    \resizebox{\columnwidth}{!}{%
        \input{tables/results_1888}
    }
\end{table}

\begin{table}[!t]
    \centering
    \caption{Solution Quality of Pipelines of RTE 2848}
    \label{tab:2848}
    \resizebox{\columnwidth}{!}{%
        \input{tables/results_2848}
    }
\end{table}

\begin{table}[!t]
    \centering
    \caption{Solution Quality of Pipelines of PEGASE 2869}
    \label{tab:2869}
    \resizebox{\columnwidth}{!}{%
        \input{tables/results_2869}
    }
\end{table}

\begin{table}[!t]
    \centering
    \caption{Solution Quality of Pipelines of RTE 6468}
    \label{tab:6468}
    \resizebox{\columnwidth}{!}{%
        \input{tables/results_6468}
    }
\end{table}

\begin{table}[!t]
    \centering
    \caption{Solution Quality of Pipelines of RTE 6470}
    \label{tab:6470}
    \resizebox{\columnwidth}{!}{%
        \input{tables/results_6470}
    }
\end{table}

\begin{table}[!t]
    \centering
    \caption{Solution Quality of Pipelines of RTE 6495}
    \label{tab:6495}
    \resizebox{\columnwidth}{!}{%
        \input{tables/results_6495}
    }
\end{table}

\begin{table}[!t]
    \centering
    \caption{Solution Quality of Pipelines of RTE 6515}
    \label{tab:6515}
    \resizebox{\columnwidth}{!}{%
        \input{tables/results_6515}
    }
\end{table}

\begin{table}[!t]
    \centering
    \caption{Solution Quality of Pipelines of PEGASE 8387}
    \label{tab:8387}
    \resizebox{\columnwidth}{!}{%
        \input{tables/results_8387}
    }
\end{table}

\begin{table}[!t]
    \centering
    \caption{Solution Quality of Pipelines of PEGASE 9241}
    \label{tab:9241}
    \resizebox{\columnwidth}{!}{%
        \input{tables/results_9241}
    }
\end{table}

\begin{table}[!t]
    \centering
    \caption{Solution Quality of Pipelines of PEGASE 13659}
    \label{tab:13659}
    \resizebox{\columnwidth}{!}{%
        \input{tables/results_13659}
    }
\end{table}

\subsection{Result Analysis}
\label{subsec:result_analysis}
Tables \ref{tab:118} to \ref{tab:13659} present numerical results for each of the cases listed in Table \ref{tab:devices}.
Each table reports: the value of $\lambda_{p}$ used in the experiment, the range of transformer tap setpoints ($\mathcal{T}$) and number of capacitor bank modules that can be activated ($\mathcal{B}$), voltage and reactive generation quality metrics (MAE$_{v}$ and MAE$_{q}$) as defined in \eqref{eq:vmae} and \eqref{eq:qmae}, the runtimes of the relaxed ($T_{r}$) and fixed ($T_{f}$) steps, and the deviation from active generation setpoint ($\Delta \PG$) and relative change in generation cost ($\%\Delta c$).
The first row of each table reports the solution quality metrics for the reference setpoint which, in this paper, corresponds to an AC-OPF solution with no adjustment of tap ratios nor capacitor banks. 
Recall that $\lamv \,{=} \,\lamq \,{=} \,\lamc \,{=} \,1$ in all experiments, while $\lamp$ varies to examine the sensitivity of voltage/VAR control to active power redispatch.

The experiments reveal clear sensitivity patterns to the $\lamp$ parameter across different test cases and device ranges. When $\lamp=\infty$ (active power dispatch fixed), voltage control effectiveness varies across systems—sometimes showing notable improvements or small gains, and occasionally leading to infeasible solutions in larger networks. By contrast, reactive power generation consistently decreases across all successful cases regardless of system scale.

As $\lamp$ decreases, voltage and VAR control performance generally improves at the cost of increased active power redispatch. Importantly, across many cases the overall generation cost decreases in a non-trivial manner (which exceeds 0.25\% on \rteXS, \rteXL, \, and \pegaseM), indicating more efficient active power dispatch. A likely driver of these cost reductions is a concurrent reduction in power losses, which are calculated but omitted from tables due to space.

Runtime scalability follows predictable trends with respect to system size and device range. The fixed stage of Algorithm \ref{alg:vvo_round_resolve} dominates total computation time, with local device ranges yielding practical solution times across all test cases. The full device range $\TAP = \pm16$ and up to 3 capacitor modules) expands the combinatorial search space, leading to longer runtimes and occasional convergence difficulties in larger systems.

A further analysis of active power losses exhibits consistent patterns across experimental settings. When $\lamp=\infty$, losses remain close to baseline values, indicating minimal impact on system efficiency. As $\lamp$ decreases and redispatch increases, losses generally decrease due to more optimal generator dispatch patterns and improved device coordination, which aligns with the observed reductions in generation cost. The trade-off between tighter voltage regulation and losses becomes more pronounced in larger systems.

Several cases present unique challenges. The RTE 6470 and 6515 systems show convergence difficulties, particularly with the full device range, reflecting the inherent computational complexity of VVO. The largest PEGASE systems (8387, 9241, 13659) exhibit mixed results, with some settings performing well while others fail to converge, highlighting the sensitivity of VVO to network size and structure.

The results demonstrate that the $\lamp$ parameter provides effective control over the trade-off between voltage quality (MAE$_v$), reactive generation (MAE$_q$), and active power redispatch ($\Delta \PG$), enabling operators to improve system efficiency and generation cost based on their specific priorities and constraints.

%% file: tables/device_composition.tex
\begin{tabular}{l rrrrr}
\toprule
\textbf{Case} & \textbf{\#buses} & \textbf{\#gen} & \textbf{\#CBs} & \textbf{\#lines} & \textbf{\#tsfm} \\
\midrule
\ieeeM    &   118 &   54 &   11 &  186 &    11 \\
% \ieeeL    & 300 & 69 & 29 & 411 & 129 \\
\pegaseXS &  1354 &  260 & 1082 &  1991 &  240 \\
\rteXS    &  1888 &  290 &   45 &  2531 &  555 \\
\rteS     &  2848 &  511 &   48 &  3776 &  783 \\
\pegaseS  &  2869 &  510 & 2197 &  4582 &  531 \\
\rteM     &  6468 &  399 &   97 &  9000 & 1580 \\
\rteL     &  6470 &  761 &   73 &  9005 & 1579 \\
\rteXL    &  6495 &  680 &   99 &  9019 & 1603 \\
\rteXXL   &  6515 &  684 &  102 &  9037 & 1615 \\
\pegaseM  &  8387 & 1865 &  487 & 14561 & 2050 \\
\pegaseL  &  9241 & 1445 & 7327 & 16049 & 2252 \\
\pegaseXL & 13659 & 4092 & 8754 & 20467 & 6675 \\
\bottomrule
\end{tabular}

%% file: tables/results_118.tex
\begin{tabular}{lll ll rr rr}
\toprule
$\lamp$ & $\TAPEDGES$ & $\CBNODES$
& MAE$_{v}$ & MAE$_{q}$
& $T_{\mathrm{r}}$ & $T_{\mathrm{f}}$
& $\Delta \PG$ & $\%\Delta\cost$ \\
\midrule
-- & $\pm 0$ & 1-1 & 0.034 & 38.00 & -- & -- & 0.00 & 0.00 \\
\midrule
% ---------- pg fixed = N ----------
$1$ & $\pm 3$ & 0-2 & 0.036 & 36.00 & 0.2 & 0.1 & 0.25 & -0.06 \\
 & $\pm 3$ & 0-3 & \textbf{0.036} & 35.51 & 0.2 & 0.1 & 0.25 & -0.06 \\
 & $\pm 16$ & 0-3 & 0.037 & 34.09 & 0.2 & 0.1 & 0.41 & \textbf{-0.07} \\
\midrule
$5$ & $\pm 3$ & 0-2 & 0.036 & 35.99 & 0.2 & 0.1 & 0.25 & -0.06 \\
 & $\pm 3$ & 0-3 & 0.036 & 35.51 & 0.2 & 0.1 & 0.25 & -0.06 \\
 & $\pm 16$ & 0-3 & 0.037 & 34.09 & 0.2 & 0.2 & 0.40 & -0.07 \\
\midrule
% ---------- pg fixed = Y ----------
$\infty$ & $\pm 3$ & 0-2 & 0.037 & 34.74 & 0.2 & 0.1 & 0.04 & -0.05 \\
 & $\pm 3$ & 0-3 & 0.037 & 34.41 & 0.2 & 0.1 & \textbf{0.04} & -0.05 \\
 & $\pm 16$ & 0-3 & 0.037 & \textbf{33.26} & \textbf{0.2} & \textbf{0.1} & 0.05 & -0.07 \\
\bottomrule
\end{tabular}%

%% file: tables/results_1354.tex
\begin{tabular}{lll ll rr rr}
\toprule
$\lamp$ & $\TAPEDGES$ & $\CBNODES$
& MAE$_{v}$ & MAE$_{q}$
& $T_{\mathrm{r}}$ & $T_{\mathrm{f}}$
& $\Delta \PG$ & $\%\Delta\cost$ \\
\midrule
-- & $\pm 0$ & 1-1 & 0.072 & 73.51 & -- & -- & 0.00 & 0.00 \\
\midrule
% ---------- pg fixed = N ----------
$1$ & $\pm 3$ & 0-2 & 0.080 & 37.12 & 4.4 & 3.5 & 1.93 & -0.18 \\
 & $\pm 3$ & 0-3 & 0.073 & 35.66 & 3.9 & 3.4 & 2.18 & -0.16 \\
 & $\pm 16$ & 0-3 & \textbf{0.066} & 39.29 & 4.8 & 3.2 & 3.90 & -0.20 \\
\midrule
$5$ & $\pm 3$ & 0-2 & 0.080 & 36.22 & \textbf{3.9} & 3.4 & 1.53 & -0.18 \\
 & $\pm 3$ & 0-3 & 0.073 & 35.50 & 4.3 & \textbf{3.2} & 2.01 & -0.15 \\
 & $\pm 16$ & 0-3 & 0.067 & 38.91 & 4.9 & 3.5 & 3.21 & \textbf{-0.20} \\
\midrule
% ---------- pg fixed = Y ----------
$\infty$ & $\pm 3$ & 0-2 & 0.078 & \textbf{23.66} & 4.2 & 3.6 & 0.15 & -0.02 \\
 & $\pm 3$ & 0-3 & 0.071 & 26.62 & 4.6 & 3.3 & \textbf{0.03} & -0.00 \\
 & $\pm 16$ & 0-3 & 0.069 & 31.48 & 5.5 & 4.7 & 0.07 & -0.01 \\
\bottomrule
\end{tabular}%

%% file: tables/results_1888.tex
\begin{tabular}{lll ll rr rr}
\toprule
$\lamp$ & $\TAPEDGES$ & $\CBNODES$
& MAE$_{v}$ & MAE$_{q}$
& $T_{\mathrm{r}}$ & $T_{\mathrm{f}}$
& $\Delta \PG$ & $\%\Delta\cost$ \\
\midrule
-- & $\pm 0$ & 1-1 & 0.040 & 30.02 & -- & -- & 0.00 & 0.00 \\
\midrule
% ---------- pg fixed = N ----------
$1$ & $\pm 3$ & 0-2 & 0.071 & 24.04 & 9.8 & 5.4 & 28.47 & \textbf{-2.19} \\
 & $\pm 3$ & 0-3 & 0.070 & 24.52 & 9.2 & 5.1 & 28.46 & -2.19 \\
 & $\pm 16$ & 0-3 & NA & NA & 10.7 & NA & NA & NA \\
\midrule
$5$ & $\pm 3$ & 0-2 & 0.071 & 23.47 & 9.1 & 5.4 & 27.06 & -2.18 \\
 & $\pm 3$ & 0-3 & 0.070 & 24.13 & \textbf{8.6} & \textbf{4.9} & 27.06 & -2.18 \\
 & $\pm 16$ & 0-3 & NA & NA & 10.2 & NA & NA & NA \\
\midrule
% ---------- pg fixed = Y ----------
$\infty$ & $\pm 3$ & 0-2 & 0.045 & 10.39 & 17.6 & 6.1 & \textbf{0.00} & -0.00 \\
 & $\pm 3$ & 0-3 & 0.045 & \textbf{9.75} & 34.9 & 4.9 & \textbf{0.00} & -0.00 \\
 & $\pm 16$ & 0-3 & \textbf{0.030} & 42.82 & 24.5 & 331.3 & 0.03 & 0.02 \\
\bottomrule
\end{tabular}%

%% file: tables/results_2848.tex
\begin{tabular}{lll ll rr rr}
\toprule
$\lamp$ & $\TAPEDGES$ & $\CBNODES$
& MAE$_{v}$ & MAE$_{q}$
& $T_{\mathrm{r}}$ & $T_{\mathrm{f}}$
& $\Delta \PG$ & $\%\Delta\cost$ \\
\midrule
-- & $\pm 0$ & 1-1 & 0.048 & 21.31 & -- & -- & 0.00 & 0.00 \\
\midrule
% ---------- pg fixed = N ----------
$1$ & $\pm 3$ & 0-2 & 0.061 & 19.36 & 10.9 & \textbf{7.2} & 0.18 & -0.06 \\
 & $\pm 3$ & 0-3 & 0.061 & 19.74 & \textbf{10.8} & 7.7 & 0.18 & -0.06 \\
 & $\pm 16$ & 0-3 & 0.063 & 19.11 & 14.1 & 9.1 & 0.36 & \textbf{-0.08} \\
\midrule
$5$ & $\pm 3$ & 0-2 & 0.061 & 19.36 & 10.8 & 8.8 & 0.11 & -0.06 \\
 & $\pm 3$ & 0-3 & 0.061 & 19.74 & 10.8 & 7.5 & 0.11 & -0.06 \\
 & $\pm 16$ & 0-3 & 0.063 & 19.11 & 13.5 & 8.6 & 0.23 & -0.08 \\
\midrule
% ---------- pg fixed = Y ----------
$\infty$ & $\pm 3$ & 0-2 & 0.058 & \textbf{17.56} & 11.9 & 9.8 & \textbf{0.05} & -0.01 \\
 & $\pm 3$ & 0-3 & \textbf{0.058} & 17.73 & 46.1 & 9.6 & 0.05 & -0.01 \\
 & $\pm 16$ & 0-3 & NA & NA & 56.0 & NA & NA & NA \\
\bottomrule
\end{tabular}%

%% file: tables/results_2869.tex
\begin{tabular}{lll ll rr rr}
\toprule
$\lamp$ & $\TAPEDGES$ & $\CBNODES$
& MAE$_{v}$ & MAE$_{q}$
& $T_{\mathrm{r}}$ & $T_{\mathrm{f}}$
& $\Delta \PG$ & $\%\Delta\cost$ \\
\midrule
-- & $\pm 0$ & 1-1 & 0.070 & 72.55 & -- & -- & 0.00 & 0.00 \\
\midrule
% ---------- pg fixed = N ----------
$1$ & $\pm 3$ & 0-2 & 0.078 & 32.91 & 13.6 & 12.7 & 0.87 & -0.14 \\
 & $\pm 3$ & 0-3 & 0.075 & \textbf{31.69} & 13.8 & 11.1 & 0.87 & -0.12 \\
 & $\pm 16$ & 0-3 & 0.073 & 35.55 & 13.9 & 11.9 & 2.19 & -0.16 \\
\midrule
$5$ & $\pm 3$ & 0-2 & 0.078 & 32.88 & \textbf{13.5} & 12.6 & 0.80 & -0.14 \\
 & $\pm 3$ & 0-3 & 0.075 & 31.88 & 14.0 & \textbf{10.8} & 0.82 & -0.12 \\
 & $\pm 16$ & 0-3 & 0.074 & 35.45 & 13.9 & 12.0 & 1.88 & \textbf{-0.17} \\
\midrule
% ---------- pg fixed = Y ----------
$\infty$ & $\pm 3$ & 0-2 & 0.066 & 34.20 & 14.0 & 61.8 & \textbf{0.00} & -0.00 \\
 & $\pm 3$ & 0-3 & 0.068 & 39.58 & 13.8 & 81.5 & 0.00 & -0.00 \\
 & $\pm 16$ & 0-3 & \textbf{0.065} & 41.64 & 82.3 & 102.1 & 0.00 & -0.00 \\
\bottomrule
\end{tabular}%

%% file: tables/results_6468.tex
\begin{tabular}{lll ll rr rr}
\toprule
$\lamp$ & $\TAPEDGES$ & $\CBNODES$
& MAE$_{v}$ & MAE$_{q}$
& $T_{\mathrm{r}}$ & $T_{\mathrm{f}}$
& $\Delta \PG$ & $\%\Delta\cost$ \\
\midrule
-- & $\pm 0$ & 1-1 & 0.048 & 34.30 & -- & -- & 0.00 & 0.00 \\
\midrule
% ---------- pg fixed = N ----------
$1$ & $\pm 3$ & 0-2 & 0.061 & 27.62 & \textbf{47.9} & 27.5 & 0.39 & -0.14 \\
 & $\pm 3$ & 0-3 & 0.059 & 28.25 & 49.6 & 29.4 & 0.47 & -0.13 \\
 & $\pm 16$ & 0-3 & 0.058 & 32.08 & 612.4 & 33.9 & 6.81 & 0.27 \\
\midrule
$5$ & $\pm 3$ & 0-2 & 0.061 & 27.62 & 51.0 & \textbf{26.9} & 0.30 & \textbf{-0.14} \\
 & $\pm 3$ & 0-3 & 0.059 & 28.22 & 55.0 & 29.3 & 0.42 & -0.13 \\
 & $\pm 16$ & 0-3 & \textbf{0.058} & 32.21 & 599.0 & 32.6 & 6.79 & 0.27 \\
\midrule
% ---------- pg fixed = Y ----------
$\infty$ & $\pm 3$ & 0-2 & 0.061 & \textbf{24.21} & 129.1 & 32.3 & 0.21 & -0.09 \\
 & $\pm 3$ & 0-3 & 0.059 & 24.41 & 127.9 & 28.4 & \textbf{0.20} & -0.09 \\
 & $\pm 16$ & 0-3 & NA & NA & 583.6 & NA & NA & NA \\
\bottomrule
\end{tabular}%

%% file: tables/results_6470.tex
\begin{tabular}{lll ll rr rr}
\toprule
$\lamp$ & $\TAPEDGES$ & $\CBNODES$
& MAE$_{v}$ & MAE$_{q}$
& $T_{\mathrm{r}}$ & $T_{\mathrm{f}}$
& $\Delta \PG$ & $\%\Delta\cost$ \\
\midrule
-- & $\pm 0$ & 1-1 & 0.048 & 26.40 & -- & -- & 0.00 & 0.00 \\
\midrule
% ---------- pg fixed = N ----------
$1$ & $\pm 3$ & 0-2 & 0.058 & 22.27 & \textbf{32.8} & 21.8 & 0.41 & -0.15 \\
 & $\pm 3$ & 0-3 & 0.057 & 21.85 & 42.4 & 23.6 & 0.38 & -0.14 \\
 & $\pm 16$ & 0-3 & 0.058 & 22.03 & 40.1 & 23.9 & 1.23 & \textbf{-0.19} \\
\midrule
$5$ & $\pm 3$ & 0-2 & 0.058 & 22.24 & 34.4 & 23.5 & 0.38 & -0.15 \\
 & $\pm 3$ & 0-3 & 0.057 & 21.82 & 40.8 & \textbf{21.4} & 0.35 & -0.14 \\
 & $\pm 16$ & 0-3 & 0.058 & 21.96 & 47.3 & 25.2 & 1.19 & -0.19 \\
\midrule
% ---------- pg fixed = Y ----------
$\infty$ & $\pm 3$ & 0-2 & 0.057 & 22.25 & 49.5 & 42.2 & 0.06 & -0.09 \\
 & $\pm 3$ & 0-3 & 0.055 & \textbf{21.62} & 47.5 & 48.8 & \textbf{0.05} & -0.08 \\
 & $\pm 16$ & 0-3 & \textbf{0.039} & 29.22 & 42.8 & 8630.6 & 0.08 & 0.12 \\
\bottomrule
\end{tabular}%

%% file: tables/results_6495.tex
\begin{tabular}{lll ll rr rr}
\toprule
$\lamp$ & $\TAPEDGES$ & $\CBNODES$
& MAE$_{v}$ & MAE$_{q}$
& $T_{\mathrm{r}}$ & $T_{\mathrm{f}}$
& $\Delta \PG$ & $\%\Delta\cost$ \\
\midrule
-- & $\pm 0$ & 1-1 & 0.050 & 34.87 & -- & -- & 0.00 & 0.00 \\
\midrule
% ---------- pg fixed = N ----------
$1$ & $\pm 3$ & 0-2 & 0.057 & 31.31 & 179.7 & 24.9 & 0.83 & \textbf{-0.28} \\
 & $\pm 3$ & 0-3 & \textbf{0.054} & 31.27 & 74.5 & \textbf{22.7} & 0.86 & -0.24 \\
 & $\pm 16$ & 0-3 & NA & NA & 51.7 & NA & NA & NA \\
\midrule
$5$ & $\pm 3$ & 0-2 & 0.057 & 31.27 & 72.0 & 23.3 & 0.82 & -0.28 \\
 & $\pm 3$ & 0-3 & 0.054 & 31.10 & 64.0 & 24.1 & 0.85 & -0.24 \\
 & $\pm 16$ & 0-3 & NA & NA & 52.7 & NA & NA & NA \\
\midrule
% ---------- pg fixed = Y ----------
$\infty$ & $\pm 3$ & 0-2 & 0.059 & \textbf{26.81} & 26.0 & 55.9 & 0.12 & -0.08 \\
 & $\pm 3$ & 0-3 & 0.055 & 27.99 & \textbf{24.7} & 30.0 & \textbf{0.08} & -0.05 \\
 & $\pm 16$ & 0-3 & NA & NA & 33.0 & NA & NA & NA \\
\bottomrule
\end{tabular}%

%% file: tables/results_6515.tex
\begin{tabular}{lll ll rr rr}
\toprule
$\lamp$ & $\TAPEDGES$ & $\CBNODES$
& MAE$_{v}$ & MAE$_{q}$
& $T_{\mathrm{r}}$ & $T_{\mathrm{f}}$
& $\Delta \PG$ & $\%\Delta\cost$ \\
\midrule
-- & $\pm 0$ & 1-1 & 0.047 & 40.03 & -- & -- & 0.00 & 0.00 \\
\midrule
% ---------- pg fixed = N ----------
$1$ & $\pm 3$ & 0-2 & 0.057 & 32.56 & 26.9 & 27.9 & 0.65 & \textbf{-0.17} \\
 & $\pm 3$ & 0-3 & 0.054 & 32.03 & 29.1 & \textbf{23.0} & 0.54 & -0.14 \\
 & $\pm 16$ & 0-3 & 0.042 & 31.01 & 38.8 & 3045.4 & 23.37 & 14.03 \\
\midrule
$5$ & $\pm 3$ & 0-2 & 0.057 & 32.47 & \textbf{25.3} & 25.7 & 0.58 & -0.17 \\
 & $\pm 3$ & 0-3 & 0.054 & 32.18 & 29.2 & 23.7 & 0.54 & -0.14 \\
 & $\pm 16$ & 0-3 & \textbf{0.034} & 30.83 & 38.9 & 2040.5 & 8.79 & 2.34 \\
\midrule
% ---------- pg fixed = Y ----------
$\infty$ & $\pm 3$ & 0-2 & 0.056 & 27.72 & 152.9 & 84.2 & 0.11 & -0.02 \\
 & $\pm 3$ & 0-3 & 0.054 & \textbf{27.27} & 165.5 & 109.3 & 0.09 & -0.02 \\
 & $\pm 16$ & 0-3 & 0.047 & 33.69 & 26.5 & 579.5 & \textbf{0.00} & -0.00 \\
\bottomrule
\end{tabular}%

%% file: tables/results_8387.tex
\begin{tabular}{lll ll rr rr}
\toprule
$\lamp$ & $\TAPEDGES$ & $\CBNODES$
& MAE$_{v}$ & MAE$_{q}$
& $T_{\mathrm{r}}$ & $T_{\mathrm{f}}$
& $\Delta \PG$ & $\%\Delta\cost$ \\
\midrule
-- & $\pm 0$ & 1-1 & 0.079 & 95.52 & -- & -- & 0.00 & 0.00 \\
\midrule
% ---------- pg fixed = N ----------
$1$ & $\pm 3$ & 0-2 & 0.083 & 79.80 & 89.1 & 186.9 & 6.28 & -1.16 \\
 & $\pm 3$ & 0-3 & 0.080 & \textbf{79.43} & \textbf{66.6} & 107.9 & 6.40 & -1.09 \\
 & $\pm 16$ & 0-3 & \textbf{0.073} & 83.20 & 113.6 & 88.2 & 16.71 & \textbf{-3.65} \\
\midrule
$5$ & $\pm 3$ & 0-2 & 0.083 & 79.88 & 84.5 & 72.7 & \textbf{4.94} & -1.15 \\
 & $\pm 3$ & 0-3 & 0.080 & 79.52 & 77.7 & 118.0 & 4.94 & -1.09 \\
 & $\pm 16$ & 0-3 & 0.073 & 82.76 & 119.0 & \textbf{71.8} & 12.68 & -3.60 \\
\midrule
% ---------- pg fixed = Y ----------
$\infty$ & $\pm 3$ & 0-2 & NA & NA & NA & NA & NA & NA \\
 & $\pm 3$ & 0-3 & NA & NA & NA & NA & NA & NA \\
 & $\pm 16$ & 0-3 & NA & NA & NA & NA & NA & NA \\
\bottomrule
\end{tabular}%

%% file: tables/results_9241.tex
\begin{tabular}{lll ll rr rr}
\toprule
$\lamp$ & $\TAPEDGES$ & $\CBNODES$
& MAE$_{v}$ & MAE$_{q}$
& $T_{\mathrm{r}}$ & $T_{\mathrm{f}}$
& $\Delta \PG$ & $\%\Delta\cost$ \\
\midrule
-- & $\pm 0$ & 1-1 & 0.060 & 62.99 & -- & -- & 0.00 & 0.00 \\
\midrule
% ---------- pg fixed = N ----------
$1$ & $\pm 3$ & 0-2 & 0.068 & \textbf{31.54} & 545.8 & \textbf{90.7} & 1.59 & -0.12 \\
 & $\pm 3$ & 0-3 & 0.063 & 32.02 & \textbf{527.0} & 117.7 & 1.47 & -0.07 \\
 & $\pm 16$ & 0-3 & 0.065 & 33.05 & 2469.6 & 160.6 & 2.05 & -0.10 \\
\midrule
$5$ & $\pm 3$ & 0-2 & 0.068 & 31.77 & 545.4 & 94.7 & 1.39 & \textbf{-0.12} \\
 & $\pm 3$ & 0-3 & 0.062 & 32.21 & 549.9 & 190.2 & 1.22 & -0.07 \\
 & $\pm 16$ & 0-3 & NA & NA & NA & NA & NA & NA \\
\midrule
% ---------- pg fixed = Y ----------
$\infty$ & $\pm 3$ & 0-2 & 0.061 & 45.16 & 867.0 & 6081.0 & \textbf{0.00} & -0.00 \\
 & $\pm 3$ & 0-3 & \textbf{0.059} & 41.75 & 1454.7 & 2719.8 & 0.00 & -0.00 \\
 & $\pm 16$ & 0-3 & 0.062 & 42.81 & 1903.1 & 1090.1 & 0.00 & -0.00 \\
\bottomrule
\end{tabular}%

%% file: tables/results_13659.tex
\begin{tabular}{lll ll rr rr}
\toprule
$\lamp$ & $\TAPEDGES$ & $\CBNODES$
& MAE$_{v}$ & MAE$_{q}$
& $T_{\mathrm{r}}$ & $T_{\mathrm{f}}$
& $\Delta \PG$ & $\%\Delta\cost$ \\
\midrule
-- & $\pm 0$ & 1-1 & 0.068 & 28.83 & -- & -- & 0.00 & 0.00 \\
\midrule
% ---------- pg fixed = N ----------
$1$ & $\pm 3$ & 0-2 & 0.074 & 15.82 & \textbf{710.5} & 109.2 & 0.61 & -0.08 \\
 & $\pm 3$ & 0-3 & 0.067 & 15.26 & 858.0 & 97.8 & 0.66 & -0.05 \\
 & $\pm 16$ & 0-3 & 0.063 & \textbf{14.18} & 2407.6 & 98.5 & 0.72 & -0.08 \\
\midrule
$5$ & $\pm 3$ & 0-2 & 0.073 & 15.92 & 774.8 & \textbf{93.7} & 0.39 & \textbf{-0.08} \\
 & $\pm 3$ & 0-3 & 0.067 & 15.34 & 786.0 & 116.3 & 0.46 & -0.05 \\
 & $\pm 16$ & 0-3 & 0.063 & 14.27 & 2426.0 & 96.9 & 0.56 & -0.08 \\
\midrule
% ---------- pg fixed = Y ----------
$\infty$ & $\pm 3$ & 0-2 & NA & NA & 1459.6 & NA & NA & NA \\
 & $\pm 3$ & 0-3 & NA & NA & 732.9 & NA & NA & NA \\
 & $\pm 16$ & 0-3 & \textbf{0.043} & 18.77 & 3742.2 & 984.1 & \textbf{0.02} & 0.02 \\
\bottomrule
\end{tabular}%

%% file: tex/conclusion.tex
\section{Conclusion}
\label{sec:conclusion}

This paper presents a Volt/VAR Optimization (VVO) framework for transmission grids that captures the operation of discrete-control devices such as capacitor banks and adjustable transformer taps.
The paper formulates the VVO problem as a mixed-integer nonlinear programming problem, which is found to be intractable for off-the-shelf MINLP solvers when solved directly.
To address this limitation, the paper proposes a heuristic relax–round–resolve algorithm to obtain feasible solutions within reasonable computing time.
The paper conducts numerical experiments on small- to large-scale transmission grids with up to 13,659 buses.
Key findings from the numerical experiments include: (i) significant improvements VAR generation across all test cases, (ii) adjusting capacitor bank and transformer taps can yield non-trivial cost savings, (iii) computational performance scales from sub-minute for smaller systems to few-minutes runtimes for large-scale cases.

These findings motivate future work to further investigate the potential of discrete-control devices in improving grid stability and reducing cost.
Future work will explore extensions to FACTS devices such as phase-shifting transformers, incorporation of device operating costs, and temporal constraints to predict control decisions for multiperiod planning. 